\documentclass[leqno,11pt]{amsart}
\usepackage{graphicx}
\usepackage{amscd}
\usepackage{amsmath}
\usepackage{caption}
\usepackage{amsfonts}
\usepackage{amssymb}
\usepackage{mathrsfs}
\usepackage{multicol}
\usepackage{color}
\usepackage{mathdots}
\usepackage{bigints}
\usepackage[labelsep=newline,font={small,up},margin=0pc]{caption}
\newcommand{\R}{\mathbb{R}}
\allowdisplaybreaks

\newcommand{\C}{{\mathbb C}}
\newcommand{\N}{{\mathbb N}}

\definecolor{blu}{rgb}{0,0,1}

\usepackage{graphicx}
\usepackage{amscd}
\usepackage{amsmath}
\usepackage{caption}
\usepackage{amsfonts}
\usepackage{amssymb}
\usepackage{mathrsfs}
\usepackage{multicol}
\usepackage{caption}
\usepackage{subcaption}
\usepackage{color}
\numberwithin{equation}{section}

\newtheorem{theorem}{Theorem}

\allowdisplaybreaks
\definecolor{darkgreen}{rgb}{0.09, 0.45, 0.27}
\definecolor{debianred}{rgb}{0.84, 0.04, 0.33}
\definecolor{orange}{rgb}{0.98, 0.6, 0.01}
\allowdisplaybreaks

\newcommand{\CC}{\mathbb C}
\newcommand{\be}{\begin{equation}}
\newcommand{\ee}{\end{equation}}
\newcommand{\bs}{\boldsymbol}

\newtheorem{lem}{Lemma}

\newtheorem{prop}{Proposition}

\allowdisplaybreaks
\definecolor{darkgreen}{rgb}{0.09, 0.45, 0.27}
\definecolor{debianred}{rgb}{0.84, 0.04, 0.33}
\allowdisplaybreaks[1]

\begin{document}

\title[On the cubic-quintic Schr\"odinger equation]{On the cubic-quintic Schr\"odinger equation}

\author[L. Baldelli]{Laura Baldelli} \email{labaldelli@ugr.es}

\address[Baldelli]{IMAG, Departamento de Análisis Matemático -- Universidad de Granada -- Campus Fuentenueva -- 18071 Granada, Spain}


\begin{abstract}
This paper explores the cubic-quintic Schr\"odinger equation in the entire Euclidean space. Our objectives are twofold: first, to advance the understanding of unresolved issues related to this equation, which are well known in the extensively studied Gross-Pitaevskii equation. Second, to consolidate existing results on the cubic-quintic equation, providing partial contributions. Specifically, we determine the explicit constant for the $L^\infty$ a priori bound and establish a partial existence result for finite energy traveling waves in suitable approximate domains of $\R^d$.  
\end{abstract}


\maketitle

\section{Introduction}

In the present paper, we consider the cubic–quintic Schr\"odinger equation 
\begin{equation}\label{cqprob}
i\partial_t\Psi-\Delta\Psi=(-\alpha_1+\alpha_3|\Psi|^2-\alpha_5|\Psi|^4)\Psi \qquad \text{in}\,\,\R^d\times\R,
\end{equation}
where $d\ge 2$, $\Psi: \R^d\times\R \to\CC$, $\alpha_1,\alpha_3,\alpha_5>0$ satisfying
\begin{equation}\label{alphai}
\alpha_3^2-4\alpha_1\alpha_5\ge 0.
\end{equation}
This condition guarantees that the polynomial $F(s):= -\alpha_1+\alpha_3 s-\alpha_5s^2$ has two distinct real positive roots $r_0^2>r_1^2>0$. In particular, we are interested in solutions of \eqref{cqprob} satisfying the following nonvanishing boundary condition $|\Psi|\to r_0$ as $|x|\to +\infty$.
Equation \eqref{cqprob} derives from the Nonlinear Schr\"odinger Equation
\begin{equation}\label{nlsprob}
i\partial_t\Psi-\Delta\Psi=F(|\Psi|^2)\Psi \qquad \text{in}\,\,\R^d\times\R
\end{equation}
deeply studied in the literature.
If $F'(r_0^2) <(>) 0$, it means that \eqref{nlsprob} is defocusing (focusing) and the sound velocity at infinity associated with \eqref{nlsprob} is $v_s=r_0\sqrt{-2F'(r_0^2)}$.
Note that the cubic nonlinearity is thereby focusing while the quintic one is defocusing, as the cubic-quintic and the Gross-Pitaevskii equation. This latter equation can be seen as \eqref{nlsprob} with a Ginzburg-Landau potential $F(s)=1-s$, namely
\begin{equation}\label{gpprob}
i\partial_t\Psi-\Delta\Psi=(1-|\Psi|^2)\Psi \qquad \text{in}\,\,\R^d\times\R,
\end{equation}
and it is one of the most studied in the literature.

At least formally, equation \eqref{nlsprob} presents two invariants: the energy and the momentum.
The energy associated with \eqref{cqprob} is 
$$E(\Psi)=\frac{1}{2}\int_{\R^d}|\nabla \Psi|^2 dx+\frac{\alpha_1}{2}\int_{\R^d}|\Psi|^2 dx-\frac{\alpha_3}{4}\int_{\R^d}|\Psi|^4 dx+\frac{\alpha_5}{6}\int_{\R^d}|\Psi|^6 dx,$$
while the momentum, which describes the evolution of the center of mass of $\Psi$, assuming that $|\Psi|\to r_0$ at infinity in a suitable sense, it is described as $\bs P =(P_1,\dots, P_d)$ with
$$P_k=\int_{\R^d} \langle i\partial_{x_k} \Psi, \Psi-r_0 \rangle,$$
where $\langle f, g\rangle:=\text{Re}(f)\text{Re}(g)+\text{Im}(f)\text{Im}(g)$ is the scalar product in $\CC$.

The cubic-quintic equation has many interpretations in physics. In the context of a Boson gas, it describes two-body attractive and three-body repulsive interactions.
Moreover, for $d=1$, this equation occurs in the description of defectons, in the theory of one-dimensional ferromagnetic and molecular chains, and in plenty of similar problems in condensed matter. On the other hand, for $d=3$, it models the evolution of a monochromatic wave complex envelope in a medium with weakly saturating nonlinearity. The propagation of stationary light beams in such a medium is also governed, provided that $d=2$ and the temporal variable $t$ is interpreted as a longitudinal coordinate. Finally, it enters in the Ginzburg-Landau two-liquid theory and equations of nuclear hydrodynamics with effective Skyrme's forces reduce quasiclassically to the same equation. For details see \cite{bp, bgmp, cs16, maris02} and the references therein.


In the last decades, particular attention has been paid to solitary waves solutions which are special classes of solutions of \eqref{nlsprob} of the form
$$\Psi(x,t) =e^{i\omega t}\psi(x_1-ct, \tilde x) \qquad \tilde x = (x_2, \dots x_d)\in\R^{d-1}.$$
Following \cite{fa03}, we can classify them as: standing waves ($c=\omega=0$), traveling waves ($\omega=0$), bound states ($c=0$).
Concerning bound states solutions to \eqref{cqprob} we refer to \cite{cs, cks, kopv17}.

In this paper, we focused on traveling waves solutions to \eqref{cqprob}, namely solutions of the form
\begin{equation}\label{ansatz}
\Psi(x,t) =\psi(x_1-ct, \tilde x) \qquad \tilde x = (x_2, \dots x_d)\in\R^{d-1},
\end{equation}
where the parameter $c\in\R$ characterizes the speed of the traveling wave, w.l.o.g. we consider $c > 0$. So that, by the ansatz \eqref{ansatz} the equation for the profile $\psi:\R^d\to \C$ in \eqref{cqprob} is
\begin{equation}\label{procq}
ic\partial_{x_1}\psi+\Delta\psi+(-\alpha_1+\alpha_3|\psi|^2-\alpha_5|\psi|^4)\psi=0.
\end{equation}
The asymptotic behavior of the Gross-Pitaevskii equation \eqref{gpprob} as $|x|\to +\infty$ was computed by Gravejat (see \cite{gra05} and references therein), confirming the earlier conjecture in \cite{jpr}. As for the asymptotic analysis of the cubic-quintic case, to the best of our knowledge, it has not yet been addressed in the literature. However, Mari\c{s} in \cite{maris13}, on page 110, noted that "it is likely that the proofs of Gravejat can be adapted to general nonlinearities".
So, following \cite{maris13} where a large class of equations are treated, by the phase invariance of the problem, we can assume that 
\begin{equation}\label{inftlim}
\psi(x)\to 1\qquad \text{as}\quad |x|\to+\infty.
\end{equation}

In particular, we are looking for traveling wave solutions with finite energy which are supposed to play an important role in the dynamics of \eqref{nlsprob}.
We will say that a solitary wave has finite energy if (and only if) the profile $\psi$ has finite energy and we will also say that a
finite-energy solitary wave is non-trivial if its profile $\psi$ is not a constant function.

At least formally, the Lagrangian associated to \eqref{procq} is defined as
$$I^c(\psi)=E(\psi)-cP(\psi),$$
where
$$E(\psi)=\frac{1}{2}\int_{\R^d}|\nabla \psi|^2 dx+\frac{\alpha_1}{2}\int_{\R^d}|\psi|^2 dx-\frac{\alpha_3}{4}\int_{\R^d}|\psi|^4 dx+\frac{\alpha_5}{6}\int_{\R^d}|\psi|^6 dx,$$
and $P:=P_1$ is the momentum with respect to the $ x_1$ variable and, under suitable integrability conditions can be written as
$$P(\psi)=-\int_{\R^d} \partial_{x_1}(Im(\psi))(Re(\psi)-1).$$
It is clear that $P$ is a well-defined functional in $H^1(\R^d)$. For a proper definition of the momentum in the larger space taken into account by Mari\c{s}, we quote \cite[Section 2]{maris13}.

This paper is dedicated to the cubic-quintic Schr\"odinger equation which is less studied in literature than the Gross-Pitaevskii one. In particular, we will mainly focus on uniform $L^\infty$ a priori estimates of any profile solution $\psi$ to \eqref{fapro}. Even if these estimates were already proved in \cite{maris08}, as the first step to prove the regularity of traveling waves, the explicit value of the constant seems to be missing: its achievement is our first aim. Through careful analysis, we can determine the explicit constant, in the sense of \cite{fa03, maris08}, by selecting the largest among the three distinct values, for details, see Section \ref{Linfty}. The exact value of the $L^\infty$ a priori bound may be used to investigate open problems related to the cubic-quintic equation.
Next, we establish the foundations to prove an existence result for finite energy traveling waves to the cubic-quintic equation in suitable approximate domains of $\mathbb{R}^d$, drawing inspiration from \cite{br}. In particular, the mountain pass geometry of the problem allows us to obtain a Palais-Smale sequence at the mountain pass level. However, it remains uncertain whether this sequence is bounded. Various techniques in the literature address this issue, such as the Pohozaev identity used in \cite{maris13} and the monotonicity trick in \cite{br}. The latter is not applicable in our case due to the lack of a definite sign of a term in the energy functional, while Mari\c{s}'s approach is not particularly relevant since it already encompasses the cubic-quintic equation. Despite numerous attempts, we have been unable to establish the boundedness of such mountain pass Palais-Smale sequence. This limitation prevents us from reaching the same conclusion as Bellazzini and Ruiz in \cite{br}, although we can demonstrate the splitting property through careful analysis and thorough work. Nevertheless, we choose to include all the aforementioned ingredients in Section \ref{existence}, providing a more comprehensive description.

Before delving into the details, the following subsections are devoted to a thorough review of the existing literature on the cubic-quintic equation, its properties such as decay and non-vanishing and to a suitable reformulation of the problem by a proper change of variable.

\subsection{Traveling waves in the literature}

Traveling wave solutions to the Schr\"odinger equation have been extensively studied in numerous papers, including \cite{bm, jpr, gr}. In particular, based on formal computations and numerical experiments, a set of conjectures, often called the Roberts program (see \cite{jpr}), has been proposed regarding the existence, qualitative properties, and stability of traveling waves in the Gross-Pitaevskii equation, and more broadly, their role in the dynamics of \eqref{nlsprob}.
\smallskip 

Concerning the existence of finite energy traveling waves, it is expected in the sub-sonic case, namely if $c\in(0, \sqrt{2})$. Many attempts have been made to prove those conjectures, see \cite{bs99, chiron04, granon, maris08}.
Finally, Mari\c{s} proved in \cite{maris13} the existence result in the entire subsonic case in dimension $d\ge 3$ for \eqref{nlsprob} with a general nonlinearity including the Gross-Pitaevskii and the cubic-quintic, by using a Pohozaev type constraint. In particular, Mari\c{s} in \cite{maris13} in the list of assumptions on the nonlinearity $F$, requires $F'(1)=-1$. So, at first sign, it seems that the cubic-quintic nonlinearity is not contained in the class considered. However, in the defocusing case, the following proper scaling, see \cite{maris13}, p. 108,
$$\Psi(x,t)=r_0\tilde\Psi\left(r_0\sqrt{-F'(r_0^2)}x,-r_0^2F'(r_0^2) t\right)$$
transform \eqref{nlsprob} in
$$i\partial_t\Psi-\Delta\Psi=\tilde F(|\tilde\Psi|)\tilde\Psi,$$
where $\tilde F(s)=-F(r_0^2s)/(r_0^2F'(r_0^2))$. So that $\tilde F'(1)=-1$. For this reason, it is possible to assume $r_0 = 1$, $F'(r_0^2)= -1$ and the speed limit as $\sqrt{2}$.

Likely ten years later, Bellazzini and Ruiz in \cite{br} confirm this result for the Gross-Pitaevskii in $d=3$ (Theorem 1.2), extend it in $d=2$ for a.e. $c\in(0,\sqrt{2})$ (Theorem 1.1) and for any $c\in(0,\sqrt{2})$ but under some assumptions on the vortex set
of the solutions (Theorem 1.3). Their approach is based on approximating domains and min-max arguments. Recently \cite{bbm} proposed another approach based on a Sobolev-type inequality which permits the confirmation of the result by Mari\c{s} in dimensions $d\ge 4$ for a large class of equations, including the Gross-Pitaevskii and the cubic-quintic. 

Unlike \eqref{gpprob}, the cubic-quintic equation exhibits unstable stationary solutions in any dimension $d \ge 1$. In \cite{lwz}, the authors use the hydrodynamic formulation to demonstrate the existence of traveling waves for cubic-quintic type equations with small speed in any dimension $d \ge 2$. This provides a simplified proof of earlier results regarding the existence of slow traveling waves, specifically in Mari\c{s}'s work \cite{maris02} for $d \ge 4$, which relied on a local version of the classical saddle-point theorem, and in an unpublished manuscript by Lin \cite{lin} for $d=2,3$. For $d=1$, explicit soliton solutions to \eqref{procq} are given in \cite{bm}, where the reader can also find details on different solution types in dimensions $d=2,3$. In the one-dimensional case, which frequently arises in Nonlinear Optics, the existence of solutions is discussed in \cite{ch12}, while the stability analysis can be found in \cite{b96}. For the two-dimensional case, traveling waves of the nonlinear Schr\"odinger equation with a general nonlinearity, including the cubic-quintic, are studied numerically in \cite{cs16}. Lastly, for $d=3$, Killip et al. in \cite{kopv12} establish the global well-posedness of the cubic-quintic nonlinear Schr\"odinger equation with non-vanishing boundary conditions at spatial infinity.

We also mention \cite{cm17} where the existence of nontrivial finite energy traveling waves for a large class of nonlinear Schr\"odinger equations of the same type considered in \cite{maris13} in space dimension $d\ge 2$ is investigated by minimizing the energy at fixed momentum (giving a set of orbitally stable traveling waves) and minimizing the action at constant kinetic energy.
Finally, in \cite{LR} the authors study positive solutions (‘ground states’) to some semilinear elliptic equations
$$-\Delta u=g(u), \qquad\text{in}\,\, \R^d\quad d\ge2$$
with a double-power nonlinearity
$$g(u)=-u^p+u^q-\mu u,\qquad p>q>1,\quad\mu>0,\quad d\ge2$$
which generalizes the cubic-quintic one which occurs if $p=5$, $q=3$.

\smallskip

Decay estimates for traveling waves in the Gross-Pitaevskii equation were established by Gravejat in \cite{gra04} (see also Remark 2.4 in \cite{br}). However, it appears that such decay properties have not been rigorously confirmed for other nonlinearities, as noted on page 232 of \cite{cm17}, where the authors addressed integrability properties rather than relying on decay estimates in order to prove a specific identity. Additionally, Remark 2.7 in \cite{lwz} suggests that the arguments presented in \cite{gra04} can be extended to prove decay for general nonlinear terms. Although we are confident that decay can be proven for generic nonlinearities, the nature of the kernel of the Gross-Pitaevskii equation appears to be insufficient for a straightforward proof for the cubic-quintic equation. Furthermore, Mari\c{s}'s approach in \cite{maris13} does not utilize decay estimates.
Moreover, the absence of decay estimates would hinder our ability to transition to the limit as $N\to+\infty$ in the event that we were able to identify a solution to problem \eqref{proN} for nearly all subsonic values of $c$. This approach, as employed by Bellazzini and Ruiz in \cite{br}, would enable the existence of a solution to \eqref{fapro} for all subsonic speeds.

\smallskip

The nonvanishing property is already proved for the Gross--Pitaevskii case in the entire Euclidean space in Proposition $2.4$ \cite{bgs9} and in slabs $\Omega_N$ in Proposition 5.1 \cite{br}. Due to the lack of decay estimates, the only setting we could consider to get the nonvanishing property is $\Omega_N$. 
However, both the approaches outlined in \cite{br} and \cite{bgs9} rely on the specific nature of the Gross-Pitaevskii operator. In particular, the Gross-Pitaevskii term in the related equation and energy always presents $(1-|\psi|^2)$, with exponent $1$ or $2$, and this fact seems to be crucial in order to conclude. While for the cubic-quintic nonlinearity, different terms appear, namely $(|\psi|^2-A)$ and $(3|\psi|^2-2A-1)$, where $A\in(0,1)$ will be specified later, and a proper comparison leads to a lower bound for $|\psi|^2$ that ultimately prevents us from reaching a conclusion. As far as we know, the non-vanishing property still remains an open problem.

\smallskip

\begin{figure}[htbp]
\centering
\includegraphics[width=8cm,  scale=0.90]{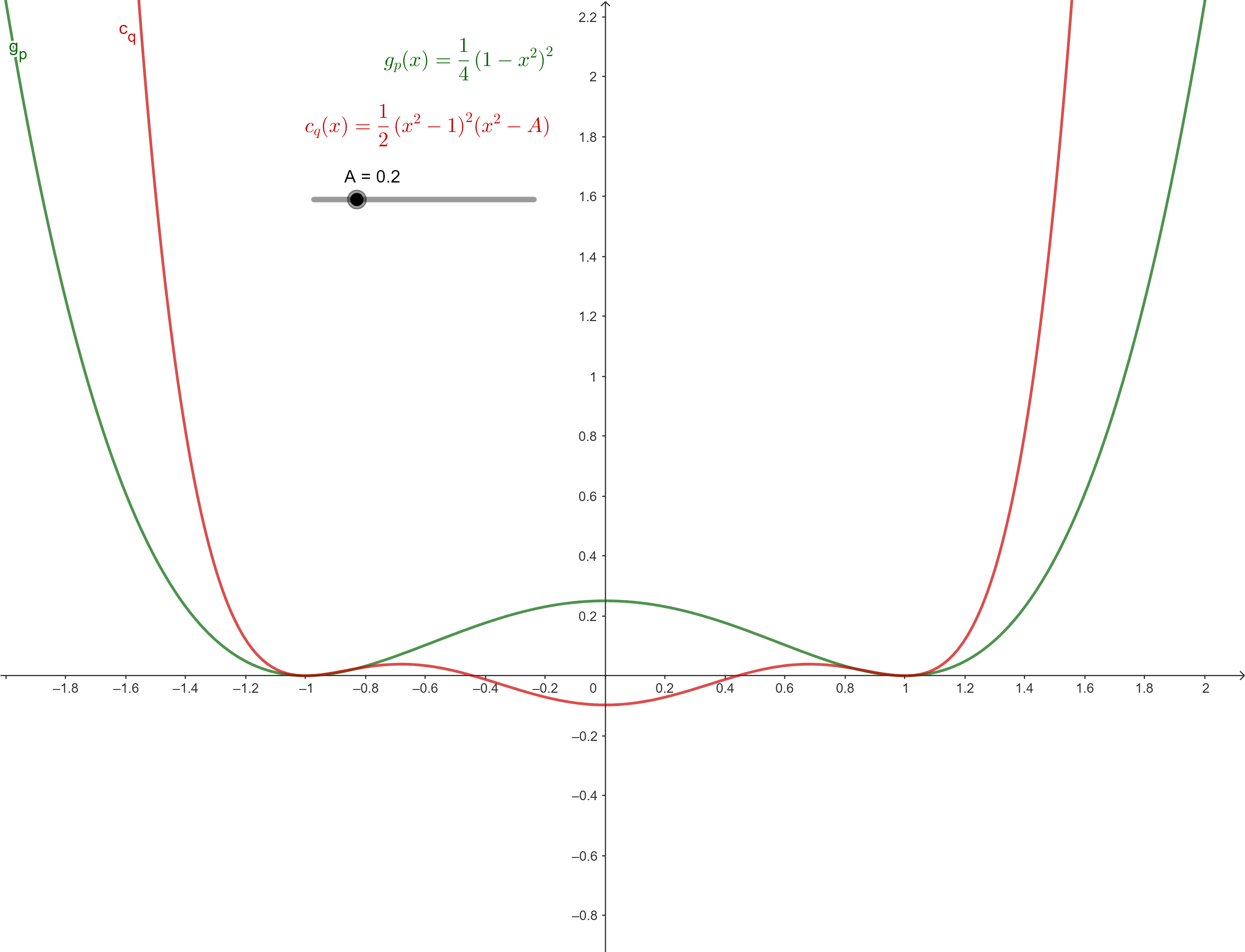}
\caption{Gross-Pitaevskii and Cubic-Quintic potentials}  \label{gpcq}
\end{figure}

Finally, by Figure \ref{gpcq}, we point out the differences and similarities between the Gross-Pitaeskii (green) and the cubic-quintic (red) potentials ($A=1/5$ in Figure \ref{gpcq}). In particular, $\psi=\pm 1$ are local minimum for the cubic-quintic potential but, while the Gross-Pitaeskii one is always nonnegative with $0$ as a local maximum, for the cubic-quintic potential $0$ is a local minimum and there exist other two local maximum between $0$ and $1$.

\subsection{A reformulation of the problem}
In what follows we will deal with an equivalent, but more useful, reformulation of \eqref{cqprob} and \eqref{procq}. 

In particular, following \cite{bp}, consider the substitution 
$$\Psi(x,t)=\gamma\Phi\left(\sqrt{\frac{\alpha_5}{2}}\gamma^2 x, \frac{\alpha_5}{2}\gamma^4 t\right), \qquad \gamma^2=\frac{3\alpha_3}{2\alpha_5(A+2r_0^2)}$$
where 
$$\frac{A}{r_0^2}=-2+\frac{3\alpha_3^2}{4\alpha_1\alpha_5}\left(1-\sqrt{1-4\frac{\alpha_1\alpha_5}{\alpha_3^2}}\right),$$
then, \eqref{cqprob} takes the form
\begin{equation}\label{facq}
i\partial_t\Phi-\Delta\Phi=\Phi(|\Phi|^2-r_0^2)(2A+r_0^2-3|\Phi|^2), \qquad \text{in}\,\,\R^d\times\R,
\end{equation}
with the boundary condition
$$\lim_{|x|\to\pm \infty}|\Phi(x,t)|=r_0.$$
So that the energy associated with \eqref{facq} becomes 
$$E(\Psi)=\frac{1}{2}\int_{\R^d}|\nabla \Psi|^2 dx+\frac{1}{2}\int_{\R^d}(|\Psi|^2-r_0^2)^2(|\Psi|^2-A) dx$$
Note that \eqref{alphai} implies $A\in(0, r_0^2)$, indeed assuming this interval for $A$, the potential term in the energy has two roots $0<A<r_0^2$.
Furthermore, by the ansatz \eqref{ansatz}, the equation for the profile $\psi$ in \eqref{facq} is given by
$$ic\partial_{x_1}\psi+\Delta\psi+\psi(|\psi|^2-r_0^2)(2A+r_0^2-3|\psi|^2)=0.$$
From now on, in order to satisfy \eqref{inftlim}, we will assume $r_0=1$. Thus the equation above becomes 
\begin{equation}\label{fapro}
ic\partial_{x_1}\psi+\Delta\psi+\psi(|\psi|^2-1)(2A+1-3|\psi|^2)=0
\end{equation}
with the related energy functional
$$
E(\psi)=\frac{1}{2}\int_{\R^d}|\nabla \psi|^2 dx+\frac{1}{2}\int_{\R^d}(|\psi|^2-1)^2(|\psi|^2-A) dx
$$
and the Lagrangian
\begin{equation}\label{eq:ac}\begin{aligned}
I^c(\psi)=E(\psi)-cP(\psi)=\frac{1}{2}\int_{\R^d}|\nabla \psi|^2 dx+\frac{1}{2}\int_{\R^d}&(|\psi|^2-1)^2(|\psi|^2-A) dx\\&+c\int_{\R^d} \partial_{x_1}(Im(\psi))(Re(\psi)-1)
\end{aligned}\end{equation}
Observe that, integrating by parts, we get
$$(I^c)'(\psi)(\phi)=\int_{\R^d} \langle \nabla\psi,\nabla\phi\rangle-c\langle i\partial_{x_1}\psi,\phi\rangle+(1-|\psi|^2)(1+2A-3|\psi|^2)\langle \psi,\phi\rangle$$
for any $\phi\in H_0^1(\R^d)$.
Moreover, the sound velocity at infinity associated to \eqref{facq} becomes $v_s=r_0\sqrt{-2F'(r_0^2)}=2\sqrt{1-A}$. Finally, we define $Q_\psi(\phi) = (I^c)''(\psi)[\phi, \phi]$, namely
\begin{equation}\label{defq}\begin{aligned}
Q_\psi(\phi)=\int_{\R^d} (|\nabla \phi|^2-c\langle\phi, i\partial_{x_1}\phi\rangle-(1-|\psi|^2)&(3|\psi|^2-1-2A)|\phi|^2\\&+4(3|\psi|^2-2-A)\langle \psi,\phi\rangle^2) dx.
\end{aligned}\end{equation}

Even if in what follows we will consider the new problem \eqref{facq}, for completeness, we report another change of variable handled in \cite{kopv12}. By rescaling (both spacetime and the values of $\Psi$ solution of \eqref{cqprob}), it suffices to consider the case $r_0^2=1$ and $\alpha_5 = 1$, so that \eqref{cqprob} reduces to 
\begin{equation}\label{kopvcq}
i\partial_t\Psi-\Delta\Psi=-\Psi(|\Psi|^2-1)(|\Psi|^2-r_1^2), \qquad \text{in}\,\,\R^d\times\R,
\end{equation}
with the boundary condition
$$\lim_{|x|\to\pm \infty}|\Psi(x,t)|=1.$$
The energy associated to \eqref{kopvcq} is 
$$E(\Psi)=\frac{1}{2}\int_{\R^d}|\nabla \Psi|^2 dx+\frac{\gamma}{4}\int_{\R^d}(|\Psi|^2-1)^2 dx+\frac{1}{6}\int_{\R^d}(|\Psi|^2-1)^3 dx$$
where $\gamma=1-r_1^2>0$, recalling that $r_1^2$ is one of the roots for $F$ under the assumption \eqref{alphai}.

\subsection{Structure of the paper}

The paper is structured as follows. Section \ref{prel} presents preliminary results, including a Pohozaev-type identity. Section \ref{Linfty} focuses on $L^\infty$ estimates for any solution to the cubic-quintic equation, and provides the explicit constant that bounds the $L^\infty$ norm. 
Finally, Section \ref{existence} paved the way to prove an existence result for the cubic-quintic equation in approximating domains, as in \cite{br}.


\section{Preliminaries}\label{prel}

In the present section, we will state some qualitative properties of traveling wave solutions to the cubic-quintic equation. 

Let us start with a Pohozaev-type identity which uses the invariance of the domain by dilations in the $\tilde{x}$ variable and it follows by \cite{maris08} and the fact that the dilations $(x_1, \tilde{x}) \mapsto (x_1, \lambda \tilde{x} )$ leave the domain $\Omega_N$ invariant, see \cite{br} for the Gross-Pitaevskii case.

\begin{lem} \label{lem:poho2}  Let $\psi$ be a finite energy solution of \eqref{fapro}. Then the following identity holds:
	$$(d-3) A(\psi)+(d-1)B(\psi)=0,$$
	where
	\[
	A(\psi)=\frac 12 \sum_{j=2}^d  \int  |\nabla_{x_j} \psi|^2  
	\]
	and
	\[
	B(\psi)=\frac{1}{2} \int |\partial_{x_1} \psi|^2 + \frac 12  \int \left(1-|\psi|^2\right)^2(|\psi|^2-A) -c\mathcal{P}(\psi). 
	\]
	Moreover, by the definition of the Lagrangian \eqref{eq:ac}, we conclude that 
	$$I^c(\psi)=\frac{2}{d-1}A(\psi)\geq 0.$$
	Finally, $I^c(\psi)=0$ if and only if $\psi$ is a constant function of modulus 1.
\end{lem}




In what follows we state two preliminary lemmas which will be useful to prove $L^\infty$ estimates in Section \ref{Linfty}. Let us start with the well-known Kato's inequality, see Lemma A in \cite{kato}.

\begin{lem}\label{kat}
If $u:\R^d\to \R$ such that $u, \Delta u\in L^1_{loc}(\R^d)$, then
$$\Delta(u^+)\ge \text{sgn}^+(u)\Delta u
,$$
in the sense of distributions.
\end{lem}

Secondly, we report a result by Brezis, see Lemma 2, p.273 in \cite{bre}.

\begin{lem}\label{bre}
Let $1<p<+\infty$ and assume $u\in L^p_{loc}(\R^d)$ satisfies
$$\Delta u\ge |u|^{p-1}u
$$
in the sense of distributions.
Then, $u\le 0$ a.e. in $\R^d$.
\end{lem}

\section{Explicit $L^\infty$ estimates}\label{Linfty}

This section is devoted to the regularity of solutions to \eqref{fapro} and the associated problem 
 restricted to slabs $\Omega_N:=\{(x_1,\tilde x)\in \R\times\R^{d-1}, \,-N<x_1<N\}$, $N\in \N$, namely
\begin{equation}\label{proN}\begin{cases}
ic\partial_{x_1}\psi+\Delta\psi+\psi(|\psi|^2-1)(2A+1-3|\psi|^2)=0\quad&\mbox{in }\Omega_N,\\\psi=1\quad&\mbox{on }\partial\Omega_N,
\end{cases}\end{equation}
In particular, we will focus on the uniform boundedness of them, taking inspiration from Farina in \cite{fa03, maris08}. For completeness, we report a regularity result by Mari\c{s}, Proposition $2.2$ in \cite{maris08}, which embraces a large class of nonlinearity, as the cubic-quintic case.

\begin{lem}\label{lem1}
Any finite energy solution $\psi$ of \eqref{fapro} or \eqref{proN} is of class $C^\infty(\R^d)$.
\end{lem}

The starting point in the proof of Lemma \ref{lem1}, that holds for any solutions, not only those with finite energy, are uniform $L^\infty$ estimates, which are valid also for problem \eqref{proN} (since the boundary condition is compatible with them).
Then, the conclusion follows by applying local elliptic regularity estimates. 

For the Gross-Pitaevskii equation \eqref{gpprob}, the $L^\infty$ estimates were proved by \cite{fa03} and the constant of the $L^\infty$ estimates is explicitly known as $\sqrt{1+c^2/4}$. Moreover, the author in \cite{fa03} at the end of the paper suggests that the technique used for the Gross-Pitaevskii equation could be applied to the cubic-quintic one to prove that "any standing (traveling/bound state) wave solution $\Psi\in L^5_{loc}(\mathbb R^d,\C)$ of \eqref{facq} satisfies $|\Psi|\le r_0$ a.e. on $\mathbb R^d$". Since Lemma \ref{lem1} does not provide explicit constants to the $L^\infty$ bound of $\psi$, we state this as our claim in the following theorem.

\begin{theorem}\label{linf}
Assume $d\ge 1$. Every solution $\psi$ of the cubic-quintic equation \eqref{fapro} satisfies the following $L^\infty$ bound
\begin{equation}\label{keyest}
|\psi(x)|\le C_A:= \sqrt{A+2}\frac{\sqrt{3+2\sqrt{3}}}{3}, \qquad \text{a.e. }\,\, x\in \R^d
\end{equation}
\end{theorem}

Note that the constant in \eqref{keyest} does not depend on the sound velocity $c$. Actually, the proof of Theorem \ref{linf} shows that one can choose other two possible constants depending on $c$, but both of them are bigger than $C_A$. In this sense, we can say that $C_A$ is the "sharp" constant.
By the way, we are not able to say if $C_A$ is "sharp" in the sense of Farina, see Remark 1.2 in \cite{fa03} because of the different nature between \eqref{facq} and the associated "Ginzburg-Landau" problem
\begin{equation}\label{gl}
\Delta w+w(|w|^2-1)(2A+1-3|w|^2)+\frac{c^2}{4}w=0,
\end{equation}
obtained by the change of variables $w=e^{i\frac{c}{2}x_1}\psi$.

Before proving the previous result, we must first establish the following key property, which is related to (H2) in \cite{maris08} and includes the explicit constant, applied to \eqref{gl}.

\begin{lem}\label{keylem}
There exists $\bar r>0$ such that
$$(s^2-1)(3s^2-2A-1)-\frac{c^2}{4}\ge 3(s-\bar r)^4$$
for any $s\ge \bar r$.
\end{lem}

\begin{proof}
Our claim is equivalent to find $\bar r>0$ such that $f(s)\ge 0$ for any $s\ge \bar r$, where
$$f(s):=12rs^3-s^2(2A+4+18r^2)+12r^3 s+2A+1-\frac{c^2}{4}-3r^4.$$
By using the software "Mathematica" it is possible to find $s_1(A,c,r), s_2(A,c,r), s_3(A,c,r)$ three zeros of $f$, such that $s_1\in \R$ and $s_2,s_3\in \C$ conjugated.
Since $f(s)\to+\infty$ as $s\to+\infty$, then two possibilities could occur:
\begin{enumerate}
\item[(i)] $f(0)<0$ and $s_2,s_3\in\C$;
\item[(ii)] $f(0)\ge0$ and $s_2=s_3\in\R$.
\end{enumerate}
The second possibility may holds for $r$ small, namely
$$0<r\le \sqrt[4]{\frac{1}{3}\left(2A+1-\frac{c^2}{4}\right)}$$
requiring that $2A+1-c^2/4\ge0$. But, from the nature of the problem, we can assume $r$ large so that (i) is in force, giving $s_1$ as the only real zero of $f$.

Now, imposing that $s_1(A,c,r)=r$ we obtain eight values of $r(A,c)$, four of them are negative or complex with negative real part, so that by the nature of $r$, only the following possibilities  remain, where
$$r_1=\sqrt{\frac{4+2A-\sqrt{4-8A+4A^2+3c^2}}{6}},\qquad r_2=\sqrt{\frac{4+2A+\sqrt{4-8A+4A^2+3c^2}}{6}}$$
$$r_3=\sqrt{A+2}\frac{\sqrt{3+2\sqrt{3}}}{3}, \qquad r_4=\sqrt{A+2}\frac{\sqrt{3-2\sqrt{3}}}{3}\in \C.$$
Note that $4-8A+4A^2+3c^2>0$ for any $A\in(0,1)$ and $c>0$. Since $r_1<r_2<r_3$, we may conclude by taking $\bar r:=r_3$.

\end{proof}

{\it Proof of Theorem \ref{linf}.}
Let $w$ be a solution of \eqref{gl} and denote $w_1:=Re(w)$, $w_2:=Im(w)$. By using Lemma \ref{kat} for $w_i-\bar r$, $i=1,2$, then using Lemma \ref{keylem} we get
$$\begin{aligned}
\Delta (w_i-\bar r)^+&\ge \text{sgn}^+(w_i-\bar r)\Delta (w_i-\bar r)\\
&=\text{sgn}^+(w_i-\bar r)\left[\left(-(|w|^2-1)(2A+1-3|w|^2)-\frac{c^2}{4}\right)w_i\right]\\
&\ge \text{sgn}^+(w_i-\bar r)(w_i-\bar r)^4w_i\ge \text{sgn}^+(w_i-\bar r)(w_i-\bar r)^5=[(w_i-\bar r)^+]^5.
\end{aligned}$$
By applying Lemma \ref{bre} since $w_i\in L^5_{loc}(\R^d)$, we have $w_i\le r$ a.e. in $\R^d$. Obviously $w$ and $-w$ satisfy \eqref{gl}. Repeating the same argument for $-w$, we infer that $-w_i\le \bar r$ a.e. in $\R^d$. Therefore, $|w_i|\le \bar r$ a.e. in $\R^d$, $i=1,2$, which implies \eqref{keyest} since $|\psi|=|w|$.

\qed

\section{An attempt to establish existence in approximating domains}\label{existence}



In this section we lay the groundwork to prove an existence result for finite energy solutions to the cubic-quintic equation in the slab $\Omega_N$, i.e. \eqref{proN}, taking inspiration from Proposition $4.1$ in \cite{br}.


The first step of our strategy consists in proving that $I^c_N$ has a mountain pass geometry on the affine space $1+H^1_0(\Omega_N)$. We drew inspiration from \cite{br}, but the distinct form of the functional and the different velocity at infinity $2\sqrt{1-A}$ make the situation quite delicate.
More precisely we aim to prove that
$$\gamma_N(c):=\inf_{g\in\Gamma(N)} \max_{t\in[0,1]} I_N^c(g(t))>0$$
where
$$\Gamma(N)=\{g\in C([0,1],(1+H_0^1(\Omega_N)) : g(0)=1, g(1)=\psi_0\},$$
with $\psi_0$ chosen so that $I_N^c(\psi_0)<0$.

\begin{prop}
Given any $c_0\in(0,2\sqrt{1-A})$, there exist $N_0>0$, $\psi_0\in 1+H_0^1(\Omega_{N_0})$ and $\chi(c_0)>0$ such that 
for any $N\ge N_0$, $c\in[c_0,2\sqrt{1-A})$:
\begin{itemize}
\item[a)] $I^c_N(\psi_0)<0$.
\item[b)] $0<\gamma_N(c)\le \chi(c_0)$.
\end{itemize}
\end{prop}

\begin{proof}
First, we want to prove that $\psi=1$ is a local minimum of the Lagrangian whenever $c^2<4(1-A)$. This could be done by proving that $(I_N^c)''(1)[\phi,\phi]=Q_1(\phi)>0$ for every $\phi\in C^\infty_0(\R^d)$, where $Q_\psi$ is defined in \eqref{defq}. Indeed, letting $\phi=u+iv$, we have
$$\begin{aligned}
Q_1(\phi)&=\int_{\R^d} (|\nabla \phi|^2-c\langle\phi, i\partial_{x_1}\phi\rangle+4(1-A)(\text{Re}(\phi))^2) dx\\
&=\int_{\R^d} (|\nabla u|^2+|\nabla v|^2+cuv_{x_1}-cu_{x_1}v+4(1-A)u^2) dx\\
&\ge\int_{\R^d} (|\nabla u|^2+|\nabla v|^2- |2cuv_{x_1}|+4(1-A)u^2) dx\\
&\ge \int_{\R^d} (|\nabla u|^2+|\nabla v|^2-v_{x_1}^2-c^2u^2+4(1-A)u^2) dx\\&\ge \int_{\R^d} (4(1-A)u^2-c^2u^2) dx>0
\end{aligned}$$
where we have used the inequality $2ab\le \lambda a^2+\frac{1}{\lambda}b^2$ with $\lambda=1/c$ and $c<2\sqrt{1-A}$.

As in \cite{br}, Lemma 4.4 in \cite{maris13}, gives a compactly supported function $\phi_0$ so that $I^{c_0}(1+\phi_0)< 0$ and $P(1+\phi_0)<0$. So it suffices to take sufficiently large $N_0$ such that $supp (\psi_0)\subset\Omega_N$, to obtain a) with $\psi_0=1+\phi_0$.

Finally, define $\gamma_0(t)=1+t\phi_0$, which obviously belongs to $\Gamma(N)$ for all $N\ge N_0$. Observe that
$$I^c_N(\gamma_0(t)) = E(\gamma_0(t))-ct^2P(\psi_0).$$
From a), we have $P(\psi_0)> 0$. Hence, for all $c \ge c_0$, 
$$I^c_N(\gamma_0(t))\le I^{c_0}_N(\gamma_0(t))\le \max_{t\in[0,1]}I^{c_0}_N\circ \gamma_0(t)=:\chi(c_0),$$
by definition. As a consequence, $\gamma_N(c)\le \chi(c_0)$ for all $N\ge N_0$, $c\ge c_0$.


\end{proof}

Thus the mountain pass geometry induces the existence of a Palais-Smale sequence at the level $\gamma_N$. Namely, a sequence $(\psi_n)_n$ such that
\begin{equation}\label{psmo}
I_N^c(\psi_n)=\gamma_N(c)+o(1), \ \ \ \|(I_N^c)'(\psi_n)\|_{H^{-1}_0(\Omega_N)}=o(1).
\end{equation}

In general, it is not clear if such Palais-Smale sequences are bounded or not.

In \cite{maris13} the author proves the boundedness of $(\psi_n)_n$ by the Pohozaev identity, here reported, for the cubic-quintic case, for completeness
\begin{equation}\label{poho}
\frac{d-3}{2(d-1)} \sum_{j=2}^d  \int_{\mathbb R^d}  |\nabla_{x_j} \psi|^2+\frac{1}{2} \int_{\mathbb R^d} |\partial_{x_1} \psi|^2 + \frac 12  \int_{\mathbb R^d} \left(1-|\psi|^2\right)^2(|\psi|^2-A) -c\mathcal{P}(\psi)=0.
\end{equation}
Indeed, one can easily prove by the use of \eqref{poho} that the sequence satisfying \eqref{psmo} has $\nabla_{x_j} \psi_n$ bounded in $L^2(\R^d)$ for any $j=2,\dots, d$. Then, following the procedure in \cite{maris13} one can also prove $\nabla_{x_1} \psi_n\in L^2(\R^d)$. This is not surprising since the Pohozaev constraint method of Mari\c{s} in \cite{maris13} also holds for the cubic-quintic equation. Indeed, it is not of interest to consider this latter technique.

Another way to prove the boundedness of $(\psi_n)_n$ is to embrace the use of the monotonicity trick argument, see \cite{br}. This result allows us to establish the boundedness of Palais-Smale sequences for almost every value of $c$, given that the Lagrangian can be decomposed as follows.
$$I^c(\psi) = A(\psi) - c B(\psi).$$
where $B(\psi)= P(\psi)$ is the momentum and $A(\psi) \geq  0$ for all $\psi \in 1+ H_0^1(\Omega_N)$.
Differently from the Gross-Pitaevskii case,
$$A(\psi)= E(\psi)=\frac{1}{2}\int_{\Omega_N}|\nabla \psi|^2 dx+\frac{1}{2}\int_{\Omega_N}(|\psi|^2-1)^2(|\psi|^2-A) dx$$
is not necessarily nonnegative and the sign assumption on $A$ is a key point for the monotonicity trick argument. Although it would be interesting to pursue the conclusion in the same way as Bellazzini and Ruiz did in \cite{br}, dealing with the cubic-quintic equation, it seems that different techniques with respect to the Pohozaev constraint, such as the monotonicity trick, is not applicable.


We conclude this section with the proof of the splitting property, the final ingredient needed for the existence of finite-energy traveling wave solutions, if only the boundedness of the sequence $(\psi_n)_n$ was available. For the sake of clarity, now we will focus on the relevant spatial dimensions, $d = 2,3$, as illustrated in the following preliminary result. The same approach can be extended to higher dimensions $d > 3$ with the necessary modifications.

\begin{lem}$[$Lemma $3.2$, \cite{br}$]$\label{eq:novan}
Let $d=2,3$ and $(Q_j)_j$ be the set of disjoint unitary cubes that covers $\Omega_N$. If $(\psi_n)_n$ is a bounded vanishing sequence, with $\psi_n=u_n+ i v_n$, in $1+ H^1_0(\Omega_n)$, i.e. such that
$$\sup_j \int_{Q_j} |u_n-1|^p+|v_n|^p  \rightarrow 0$$
as $n\to+\infty$, for some $2\leq p<+\infty$ if $d=2$, $2\leq p<6$ if $d=3$, then
$$\int_{\Omega_N} |u_n-1|^r+|v_n|^r  \rightarrow 0$$
as $n\to+\infty$, for any $2< r<+\infty$ if $d=2$, $2<r<6$ if $d=3.$
\end{lem}

\begin{prop}[Splitting property]
Given $0 <c<2\sqrt{1-A}$ and $(\psi_n)_n$ a bounded Palais-Smale sequence at the energy level $\gamma_N(c)$. Then there exist $k$ sequences of points $(y_n^j)_j \subset \{0\} \times \R^{d-1}$, $1\leq j\leq k,$ with $|y_n^j - y_n^k| \to +\infty$ if $j \neq k$, such that, up to subsequence,
$$\psi_n -1  =w_n + \sum_{j=1}^k (\psi^j(\cdot + y_n^j)-1)  \text{ with } w_n \rightarrow 0 \text { in }  H^1_{0}(\Omega_N),$$
\begin{equation}\label{eq:spl12}
\|\psi_n-1\|_{ H^1_0(\Omega_N)}^2\rightarrow \sum_{j=1}^k  \|\psi^j-1\|_{H^1_0(\Omega_N)}^2,
\end{equation}
$$I^c_N(\psi_n ) \rightarrow \sum_{j=1}^k I^c_N(\psi^j),$$
as $n\to+\infty$, where $\psi^j$ are nontrivial finite energy solutions to \eqref{proN}. 

In particular $I^c_N(\psi^j)\leq \gamma_N(c) \leq \chi(c_0)$, $E({\psi^j}) \leq \limsup_{n \to +\infty}  \, E(\psi_n)$ for all $j =1, \dots k$.
\end{prop}

\begin{proof}

We first claim that $\psi_n$ is not vanishing. By a contradiction argument, assume that $\psi_n$ vanishes, so Lemma \eqref{eq:novan} implies
\begin{equation}\label{eq:consvan}
\int_{\Omega_N} |u_n-1|^r+|v_n|^r  \rightarrow 0
\end{equation}
as $n\to +\infty$, for any $2< r<+\infty$ if $d=2$, $2<r<6$ if $d=3.$ We first claim that 
\begin{equation}\begin{aligned}\label{eq:crucvan}
\int_{\Omega_N} (1-u_n^2-v_n^2)^2&(u_n^2+v_n^2-A)=\int_{\Omega_N} 4(1-A)(u_n-1)^2\\&+\int_{\Omega_N}3(u_n-1)^4v_n^2+3(u_n-1)^2v_n^4+(u_n-1)^6+v_n^6+o(1).
\end{aligned}\end{equation}
as $n\to +\infty$. Indeed, notice that
\begin{equation}\label{eq1s}
(1-u_n^2-v_n^2)^2=4(u_n-1)^2+(u_n-1)^4+v_n^4+4(u_n-1)^3+4(u_n-1)v_n^2+2(u_n-1)^2v_n^2
\end{equation}
and 
$$(u_n^2+v_n^2-A)=v_n^2+(u_n-1)^2+2(u_n-1)+1-A.$$
By using H\"older's inequality and \eqref{eq:consvan}, we get
$$
\int_{\Omega_N}4(u_n-1)^2(u_n^2+v_n^2-A)=\int_{\Omega_N} 4(1-A)(u_n-1)^2+o(1)
$$
as $n\to+\infty$. By similar computations, we arrive to prove \eqref{eq:crucvan}.
Therefore, recalling the definition of the Lagrangian in \eqref{eq:ac}, we have
\begin{equation}\begin{aligned} \label{eq:vanen}
I_N^c(u_n,v_n)&=\frac{1}{2} \int_{\Omega_N} |\nabla u_n|^2 + \frac{1}{2} \int_{\Omega_N} |\nabla v_n|^2 -c \int_{\Omega_N} (1-u_n)\partial_{x_1} v_n \\
&\quad+\frac 12\int_{\Omega_N} (1-u_n^2-v_n^2)^2(u_n^2+v_n^2-A)\\&=\frac{1}{2} \int_{\Omega_N} |\nabla u_n|^2 + \frac{1}{2} \int_{\Omega_N} |\nabla v_n|^2 -c \int_{\Omega_N} (1-u_n)\partial_{x_1} v_n \\
&\quad+ \int_{\Omega_N} 2(1-A)(u_n-1)^2+\int_{\Omega_N}\frac{3}{2}(u_n-1)^4v_n^2+\frac{3}{2}(u_n-1)^2v_n^4\\
&\quad+\frac{1}{2}(u_n-1)^6+\frac{1}{2}v_n^6 +o(1). 
\end{aligned}\end{equation}
as $n\to +\infty$. On the other hand, direct computation gives
\begin{equation}\label{incp}\begin{aligned}
o(1)=(I_N^c)'&[\psi_n](1-\psi_n)=\int_{\Omega_N}|\nabla u_n|^2+|\nabla v_n|^2  -2c \int_{\Omega_N} (1-u_n)\partial_{x_1} v_n(x) \\&-\int_{\Omega_N} (1-u_n^2-v_n^2)(1+2A-3(u_n^2+v_n^2))(u_n(1-u_n)-v_n^2),\end{aligned}
\end{equation}
as $n\to +\infty$. Noting that
$$1-u_n^2-v_n^2=-(u_n-1)^2-2(u_n-1)-v_n^2,$$
\begin{equation}\label{eq2s}
1+2A-3(u_n^2+v_n^2)=-3(u_n-1)^2-6(u_n-1)-3v_n^2+2(A-1),
\end{equation}
arguing as before we notice that
$$\int_{\Omega_N} (1-u_n^2-v_n^2)(1+2A-3(u_n^2+v_n^2))v_n^2= \int_{\Omega_N}3(u_n-1)^4v_n^2+6(u_n-1)^2v_n^4+3v_n^6+o(1),$$
as $n\to +\infty$. Similarly.
$$\begin{aligned}
\int_{\Omega_N}& (1-u_n^2-v_n^2)(1+2A-3(u_n^2+v_n^2))(u_n(1-u_n))\\
&=\int_{\Omega_N} (1-u_n^2-v_n^2)(1+2A-3(u_n^2+v_n^2))(1-u_n)\\
&\qquad-\int_{\Omega_N} (1-u_n^2-v_n^2)(1+2A-3(u_n^2+v_n^2))(1-u_n)^2\\
&=-4(1-A)\|u_n-1\|_{L^2(\Omega_N)}^2-\int_{\Omega_N}3(u_n-1)^6+6(u_n-1)^4v_n^2+3(u_n-1)^2v_n^4+o(1),
\end{aligned}$$
as $n\to +\infty$. We get hence that
\begin{equation}\begin{aligned}\label{eq:vander}
(I_N^c)'&[\psi_n](1-\psi_n)=\int_{\Omega_N}|\nabla u_n|^2+|\nabla v_n|^2  -2c \int_{\Omega_N} (1-u_n)\partial_{x_1} v_n(x) \\
&\qquad+4(1-A)\|u_n-1\|_{L^2(\Omega_N)}^2+ \int_{\Omega_N}3(u_n-1)^4v_n^2+6(u_n-1)^2v_n^4+3v_n^6  \\
&\qquad\qquad+\int_{\Omega_N}6(u_n-1)^4v_n^2+3(u_n-1)^2v_n^4+3(u_n-1)^6+o(1)\\
&=\int_{\Omega_N}|\nabla u_n|^2+|\nabla v_n|^2  -2c \int_{\Omega_N} (1-u_n)\partial_{x_1} v_n(x)+4(1-A)\|u_n-1\|_{L^2(\Omega_N)}^2 \\
&\qquad+ \int_{\Omega_N}9(u_n-1)^4v_n^2+9(u_n-1)^2v_n^4+3v_n^6+3(u_n-1)^6+o(1),
\end{aligned}\end{equation}
as $n\to +\infty$. Taken into account  \eqref{eq:vanen} and\eqref{eq:vander} we conclude
$$\begin{aligned}
\gamma(N)+o(1)&=I_N^c(\psi_n)-\frac 12 (I_N^c)'[\psi_n](1-\psi_n)
\\&=- \int_{\Omega_N}3(u_n-1)^4v_n^2+3(u_n-1)^2v_n^4+v_n^6+(u_n-1)^6+o(1),
\end{aligned}$$
as $n\to +\infty$, which is a contradiction, since $\gamma(N)>0$. 

Thus, we excluded vanishing. Thus, there exists a sequence  $y^1_n \in \{0\} \times \R^{d-1}$ and $\psi^1 \in 1 + H_0^1(\Omega_N)$, $\psi^1 \neq 1$,  such that
$$\psi_n(\cdot + y^1_n) - \psi^1 \rightharpoonup 0 \mbox { in } H_0^1(\Omega_N)$$
up to a subsequence. By the invariance under translation, we can assume that $y_n^1=0$. Define $z_n^1=\psi_n-\psi^1$. From the definition of weak convergence, we obtain that
$$\begin{aligned}\frac 12 \int_{\Omega_N}|\nabla \psi_n|^2 -c \mathcal{P}(\psi_n)=\frac 12 \int_{\Omega_N}&|\nabla \psi^1|^2 -c \mathcal{P}(\psi^1)\\
&+\frac 12 \int_{\Omega_N}|\nabla z_n^1|^2 -c \mathcal{P}(1+z_n^1)+o(1).
\end{aligned}$$
as $n\to +\infty$. By H\"older's inequality and some calculations, we notice that the nonlinear term fulfills the following splitting property
$$\begin{aligned}
\int_{\Omega_N} \left(1-|\psi_n|^2\right)^2\left(|\psi_n|^2-A\right) =\int_{\Omega_N}& \left(1-|\psi^1|^2\right)^2\left(|\psi^1|^2-A\right) \\&+\int_{\Omega_N} \left(1-|1+z_n^1|^2\right)^2\left(|1+z_n^1|^2-A\right) +o(1).
\end{aligned}$$
as $n\to +\infty$. As a consequence, the norm and the action splits as 
$$\|\psi_n-1\|^2=\|\psi^1-1\|^2+\|z_n^1\|^2+o(1).$$
\begin{equation}\label{eq:split}
I_N^c(\psi_n ) = I_N^c(\psi^1)+ I_N^c(\psi_n-\psi^1)+o(1),
\end{equation}
as $n\to +\infty$.

Clearly, $\psi^1$ is a weak solution of \eqref{proN}, i.e.
\begin{equation}\label{inccp}\begin{aligned}
(I_N^c)'&[\psi^1](1-\psi^1)=\int_{\Omega_N}|\nabla u^1|^2+|\nabla v^1|^2  -2c \int_{\Omega_N} (1-u^1)\partial_{x_1} v^1(x) \\&-\int_{\Omega_N} (1-(u^1)^2-(v^1)^2)(1+2A-3((u^1)^2+(v^1)^2))(u^1(1-u^1)-(v^1)^2)=0.\end{aligned}
\end{equation}
Now if $z_n^1 \rightarrow 0$ in $H^1_0(\Omega_N)$ the lemma is proved. Let us assume the contrary, i.e. that  $z_n^1 \rightharpoonup 0$ and
$z_n^1\nrightarrow 0$ in $H^1_0(\Omega_N)$. 

We will prove that there exists a sequence of points $y^2_n\in \{0\} \times \R^{d-1}$, $|y^2_n|\to+\infty$ and $\psi^2 \in  H_0^1(\Omega_N)$, $\psi^2 \neq 1$, such that
$z_n^1(\cdot +y^2_n)+1 - \psi^2 \rightharpoonup 0$ in $H_0^1(\Omega_N)$.  Let us argue again by contradiction assuming that the sequence $z_n^1$ vanishes which implies 
by Lemma \ref{eq:novan} that
\begin{equation*}
\int_{\Omega_N} |u_n-u^1|^r+|v_n-v^1|^r  \rightarrow 0
\end{equation*}
as $n\to +\infty$, for any $2< r<+\infty$ if $d=2$, $2<r<6$ if $d=3,$
where $\psi^1 = u^1 + i v^1$.
Recalling \eqref{incp}, \eqref{inccp}, using the splitting property \eqref{eq:split} and $(I_N^c)'[\psi_n](1-\psi_n)-(I_N^c)'[\psi^1](1-\psi^1)=o(1)$ as $n\to+\infty$ we get
$$\begin{aligned}\int_{\Omega_N}&|\nabla z_n^1|^2 -2c \int_{\Omega_N} Re(z_n^1)\partial_{x_1}Im( z_n^1(x)) +o(1)\\
&=J:=\int_{\Omega_N} (1-u_n^2-v_n^2)(1+2A-3(u_n^2+v_n^2))(u_n(1-u_n)-v_n^2)\\
&\qquad\qquad-\int_{\Omega_N} (1-(u^1)^2-(v^1)^2)(1+2A-3((u^1)^2+(v^1)^2))(u^1(1-u^1)-(v^1)^2)\end{aligned}$$
Now, using the elementary inequality 
\begin{equation}\label{ineqel}
c x y\leq \frac{c^2}{4}x^2+y^2
\end{equation}
we have
\begin{equation}\label{eq:import2}
\int_{\Omega_N}|\nabla Re(z_n^1)|^2+ \left(1-\frac{c^2}{2}\right)\int_{\Omega_N}|\nabla Im(z_n^1)|^2 \leq 2 \int_{\Omega_N}|Re(z_n^1)|^2+ J+o(1),
\end{equation}
as $n\to +\infty$. Note
$$\begin{aligned}
&(1-u_n^2-v_n^2)(1+2A-3(u_n^2+v_n^2))(u_n(1-u_n)-v_n^2)\\&=(1-u_n^2-v_n^2)^2(1+2A-3(u_n^2+v_n^2))+(1-u_n^2-v_n^2)(u_n-1)(1+2A-3(u_n^2+v_n^2)).
\end{aligned}$$
Recalling \eqref{eq1s}, \eqref{eq2s} and
\begin{equation}\label{1-psipsi}
(1-u_n^2-v_n^2)(u_n-1)=-2(1-u_n)^2+(1-u_n)^3+v_n^2(1-u_n),
\end{equation}
assuming that $z_n^1=\psi_n-\psi^1\rightharpoonup 0$, we get
$$2 \int_{\Omega_N}|Re(z_n^1)|^2+J=o(1)$$
as $n\to +\infty$, and hence we get a contradiction with \eqref{eq:import2}.
We have hence proved the existence of a sequence $y_n^2 \in \{0\} \times \R^{d-1}$ and $\psi^2 \in 1 + H_0^1(\Omega_N)$, $\psi^2 \neq 1$, such that 
$$z_n^1(\cdot +y_n^2)+1- \psi^2 \rightharpoonup 0  \mbox { in } H_0^1(\Omega_N).$$
Clearly, $|y_n^2| \to +\infty$, and $\psi^2$ is also a (weak) solution of \eqref{proN}. We now define
$$z_n^2=z_n^1(\cdot)+1-\psi^2(\cdot -y_n^2).$$
Now we iterate the splitting argument to obtain
$$\|\psi_n-1\|^2=\|\psi^1-1\|^2+\|\psi^2-1\|^2+\|z_n^2\|^2+o(1),$$
\begin{equation}\label{3.16}
I_N^c(\psi_n ) = I_N^c(\psi^1)+I_N^c(\psi^2)+ I_N^c(\psi_n-\psi^1)+o(1),
\end{equation}
as $n\to +\infty$. If $z_n^2 \rightarrow 0$ in $H^1_0(\Omega_N)$ the lemma is proved. Instead, if $z_n^1 \nrightarrow 0$ we can repeat the previous procedure to find
 $y^3_n\in \{0\} \times \R^{d-1}$, $|y^3_n|\to+\infty$, $|y^3_n-y^2_n|\to+\infty$ and $\psi^2 \in  H_0^1(\Omega_N)$, $\psi^2 \neq 1$, and $z_n^3, \psi^3$ in $H^1_0(\Omega_N)$ as before.

We aim to show that we can have only a finite number of iterative steps. On this purpose, we claim that 
$$\inf_{\psi \in \mathcal{N}}\|1-\psi\|_{H^1_0(\Omega_N)}>0$$
where
$$\mathcal{N}:=\left\{ \psi \in  1 + H^1_0(\Omega_N), \psi\neq 0, \ (I_N^c)'[\psi](1-\psi)=0 \right\}.$$
This allows us to conclude thanks to \eqref{eq:spl12}. Recalling \eqref{eq1s} and \eqref{1-psipsi}, the following identity is in force
$$\begin{aligned}
  \int_{\Omega_N} &(1-u^2-v^2)(u(1-u)-v^2)(1+2A-3(u^2+v^2))=\int_{\Omega_N} \bigl[2(u-1)^2+3(u-1)^3 \\
&+   3 v^2(u-1)+2(u-1)^2v^2+ (u-1)^4+v^4\bigr](1+2A-3(u^2+v^2)),
\end{aligned}$$
Now, thanks to the inequality \eqref{ineqel}, we have
$$ 2c \int_{\Omega_N} (1-u)\partial_{x_1} v  \leq   \frac{c^2}{4(1-A)} \int_{\Omega_N} |\nabla v|^2 + 4(1-A) \int_{\Omega_N} (u-1)^2 $$
so that, by using also \eqref{eq2s}, 
$$\begin{aligned}
0=(I_N^c)'[\psi]&(1-\psi)\geq  \int_{\Omega_N} |\nabla u|^2 +\left(1-\frac{c^2}{4(1-A)}\right) \displaystyle  \int_{\Omega_N} |\nabla v|^2  \\
& \quad- 4(1-A) \int_{\Omega_N} (u-1)^2 + \int_{\Omega_N} (3v^2+3(u-1)^2+6(u-1)+2(1-A)) \\
&\qquad\cdot\bigl[2(u-1)^2+3(u-1)^3+   3 v^2(u-1)+2(u-1)^2v^2+ (u-1)^4+v^4\bigr]\\
\end{aligned}$$
The inequality above implies
$$\begin{aligned}
C_1\|\psi -1\|_{H^1_0(\Omega_N)}^3+C_2\|\psi -1\|_{H^1_0(\Omega_N)}^4+C_3&\|\psi -1\|_{H^1_0(\Omega_N)}^5+C_4\|\psi -1\|_{H^1_0(\Omega_N)}^6\\&\qquad\geq \left(1-\frac{c^2}{4(1-A)}\right)\|\psi - 1\|_{ H^1_0(\Omega_N)}^2
\end{aligned}$$
with $C_i>0$ and hence $\inf_{\psi \in \mathcal{N}}\|\psi-1\|_{ H^1_0(\Omega_N)}>0$.\\
Finally, recall that $I_N(\psi^j )>0$ by Lemma \ref{lem:poho2}, which holds also for any solutions to \eqref{proN}. Now, by \eqref{3.16}, we have
$$\gamma_N+o(1)=I_N^c(\psi_n)=\sum_j I_N^c(\psi^j(\cdot +y_n^j)) +o(1) \geq I_N^c(\psi^j) +o(1)$$
and hence $I_N^c(\psi^j)\leq \gamma_N.$

\end{proof}

Note that even if the counterpart of Proposition 4.1 in \cite{br} for \eqref{proN} was established, the absence of the non-vanishing property would prevent us from obtaining uniform energy bounds independent of the parameter $N$. Consequently, this obstructs the passage to the limit as $N \to +\infty$, which is necessary to get the existence of finite-energy traveling waves in the full subsonic range for \eqref{procq} in the entire Euclidean space $R^d$.




\section*{Acknowledgments}
The author warmly thanks Professor David Ruiz for his encouragement in undertaking this project and for his valuable comments.\\
The author is a member of the {\em Gruppo Nazionale per l'Analisi Ma\-te\-ma\-ti\-ca, la Probabilit\`a e le loro Applicazioni} (GNAMPA) of the {\em Istituto Nazionale di Alta Matematica} (INdAM). She is partially supported by the ``Maria de Maeztu'' Excellence Unit IMAG, reference CEX2020-001105-M, funded by MCIN/AEI/10.13039/501100011033/ and by the INdAM-GNAMPA Project 2024 titled {\em Regolarità ed esistenza per operatori anisotropi} (E5324001950001).

\end{document}